\documentclass[12pt]{amsart}
\usepackage{amssymb,amsmath,amsthm,amscd,mathrsfs}

\newtheorem{teo}{Theorem}[section]
\newtheorem{pro}[teo]{Proposition}
\newtheorem{lem}[teo]{Lemma}
\newtheorem{cor}[teo]{Corollary}

\theoremstyle{definition}

\theoremstyle{remark}
\newtheorem{rem}[teo]{Remark}

\begin{document}

\title[Base locus of the generalized theta divisor]{Some examples of vector bundles in the base locus of the generalized theta divisor}

\author[Casalaina-Martin]{Sebastian Casalaina-Martin }
\thanks{The first author was partially supported by NSF MSPRF grant DMS-0503228.  He would like to thank Harvard University for its hospitality during the preparation of this paper.}
\address{Sebastian Casalaina-Martin, Department of Mathematics, University of Colorado at
Boulder, Campus Box 395, Boulder, CO 80309-0395, USA}
\email{casa@math.colorado.edu}

\author[Gwena]{Tawanda Gwena}
\address{Tawanda Gwena, Mathematics Department, Tufts University, Medford MA 02155, USA}
\email{tawanda.gwena@tufts.edu}

\author[Teixidor i Bigas]{Montserrat Teixidor i Bigas}
\address{Montserrat Teixidor i Bigas, Mathematics Department, Tufts University, Medford MA 02155, USA}
\email{montserrat.teixidoribigas@tufts.edu}

\date{\today}
\begin{abstract}
We show that on the moduli space of semi-stable vector bundles of
fixed rank and determinant (of any degree)
 on a curve, the base locus of the theta divisor  as well as its $n$-multiples is large.
  This extends known results
  for  the case of trivial determinant and $n=1$.

\end{abstract}
\maketitle

\section*{Introduction}
Let $C$ be a smooth projective curve over $\mathbb C$ of genus
$g(C)=g\ge 2$.   Fix $r,d\in \mathbb Z$ and $r\ge 1$.
 Consider  the moduli  space  $\mathscr U(r,d)$ of (equivalence classes of) semi-stable vector bundles on $C$
 of rank $r$ and degree $d$.
 For a line bundle $L\in \operatorname{Pic}^d(C)$, consider the subspace $S\mathscr U(r,L)\subset \mathscr U(r,d)$
  of vector bundles with determinant $L$. Note that because of the
  isomorphisms given by tensoring with a line bundle,
there are at most $r$ not necessarily isomorphic moduli spaces
$\mathscr U(r,d)$. Therefore, it suffices to consider the slope
modulo integers. In the case of fixed determinant, only the degree
of $L$ is relevant.

There is an ample line bundle $\mathscr L$ on $S\mathscr U(r,L)$,
called the determinant line bundle, such that
$\operatorname{Pic}(S\mathscr U(r,L))=\mathbb Z \langle\mathscr L
\rangle$ \cite{DN}.  We will be interested in the base locus of the
linear system $|\mathscr L^{\otimes n}|$ on $S\mathscr U(r,L)$, for
positive integers $n$. When $d=0$, (or equivalently when $r$ divides
$d$), a great deal is known about this base locus. In \cite{R},
Raynaud provided a finite number of vector bundles in the base locus
for very high rank with respect to the genus. Later, Popa \cite{P}
showed that there are positive dimensional families of such
examples. Some improvements on those results were obtained in
\cite{A}, \cite{gt}. Finally, \cite{H} gave, at least theoretically,
a characterization of all the points in the base locus (see also
\cite{H2}).

In sharp contrast with the degree zero case, not much has been done
for other degrees, the main reason being that a characterization of
fixed points was not available untill the Strange Duality Conjecture
was proved in \cite{mo}(see also \cite{B}). We generalize these
results in two directions. We prove that, for any slope, one can
find a positive dimensional family in the base locus of the theta
divisor in suitable rank. We also show that the base loci of
$|\mathscr L^{\otimes n}|$ with $n>1$ are often quite large. Some
remarks in this direction were provided in \cite{P}.

\begin{teo}\label{teomain}
Fix a slope $\mu\in \mathbb Q$, and a level $n\in \mathbb N$.
For $N\in \mathbb N$ and $L\in \operatorname{Pic}^{\mu N}(C)$,  the
base locus of the linear system
 $|\mathscr L^{\otimes n}|$
on $S\mathscr U (N,L)$ is large provided
 $N>>0$ and $g>>0$.
\end{teo}

For more precise conditions on $N$ and $g$, and for estimates on the dimension of the base loci, see Theorem \ref{teo1} and Theorem \ref{teo2}.
We point out that Theorem 8.1 in \cite{pr} states that on $S\mathscr U(N,\mathscr O_C)$, $\mathscr L^{\otimes n}$ is
globally generated for $n\ge \lfloor\frac{N^2}{4}\rfloor$, and it has been conjectured to be  for
$n\ge N-1$.  The bound obtained in these notes  is far from this.

 The main idea
is to construct explicit examples of vector bundles as in \cite{gt},
\cite{P} by taking wedge powers
 of certain vector bundles obtained as kernels of an evaluation map.  A technical point
  (cf, Prop. \ref{promainexa}) using results in \cite{rt} and
  \cite{mo}
  shows that these are in the base locus.

\section{Preliminaries}
Given $L\in \operatorname{Pic}^d(C)$ and $F\in \mathscr U(r_F,d_F)$
such that $\mu(\xi\otimes F)=g-1$ for $\xi\in S\mathscr U(r,L)$  set
$ \Theta_F=\{\xi\in S\mathscr U (r,L):h^0(\xi\otimes F)>0\}. $ This
definition of $ \Theta_F$ depends on the choice of suitable $r,L$
that, for simplicity, are not included in the notation.

From $\mu(\xi\otimes F)=g-1$, $r_F=n\frac{r}{gcd(r,d)}$ for some
$n\in \mathbb N$. It is well known that if $\Theta_F$ is a divisor,
then $\Theta_F\in |\mathscr L^{\otimes n}|$. The Strange Duality
Theorem (Theor. 4 \cite{mo}) implies that the linear system
$|\mathscr L^{\otimes n}|$ is spanned by such divisors, and thus its
base locus equals
$$ \{\xi\in S\mathscr U(r,L):h^0(F\otimes \xi)>0 \
\forall \ F \in \mathscr U(r_F,d_F)\}.$$

\begin{lem}
Let $L\in \operatorname{Pic}^d(C)$, and suppose that on
$S\mathscr U(r,L)$
the base locus of $|\mathscr L^{\otimes n}|$ has dimension $b$.  If $r'\ge r$, $L'\in \operatorname{Pic}^{d'}(C)$,
and $\frac{d}{r}-\frac{d'}{r'}\in \mathbb Z$,
 then the base locus of $|\mathscr L|^{\otimes n}$ on $S\mathscr U(r',L')$ has dimension at least
$ b+(r'-r)^2(g-1)+1-\delta_{r,r'}$,  where $\delta_{r,r'}$ is the
Kronecker delta.
     \end{lem}
\begin{proof}
It suffices to consider slopes up to integers. If
$\frac{d}{r}=\frac{d'}{r'}$, consider direct sums of vector bundles
of slope $\frac{d}{r}$, and of the appropriate ranks (see \cite{gt}
Corollary 5.3).
\end{proof}

Given a globally generated vector bundle $E$ on $C$, define a vector
bundle $M_E$ as the kernel of the evaluation map
$$
0\rightarrow M_E \rightarrow H^0(E)\otimes \mathscr O_C \rightarrow E \rightarrow 0.
$$

From \cite{gt} Thm.~1.1 and 1.2, if $d_E-2r_Eg> -g$ and
$0<i<d_E-r_Eg$, the map  given by $E\mapsto \wedge^iM_E^\vee$ is
generically finite.  Below, we give some conditions that  insure
that for each $E\in \mathscr U(r_E,L)$, the corresponding $\wedge^i
M_E^\vee$ is in the base locus of $|\mathscr L^{\otimes n}|$ on $
S\mathscr U\left(\binom{d_E-r_Eg}{i},L^{\otimes
\binom{d_E-r_Eg-1}{i-1}}\right)$ for certain $n$.

\section{Examples}
\begin{pro}\label{promainexa}
If $E\in \mathscr U(r_E,d_E)$, and $r_F,i\in \mathbb N$ satisfy:
\begin{enumerate}
\item $\frac{r_Fid_E}{d_E-r_Eg}\in \mathbb N,$
\item $\frac{g-1}{g}\cdot\frac{d_E-r_Eg}{r_Er_F}\ge i\ge g$,
\end{enumerate}
then setting  $d_F=r_F(g-1)-\frac{r_Fid_E}{d_E-r_Eg}$,
$$\Theta_{\wedge^iM_E^\vee} :=\{F\in \mathscr U(r_F,d_F):h^0(F\otimes \wedge^iM_E^\vee)>0\}=\mathscr U(r_F,d_F).$$
     \end{pro}
\begin{proof}
Condition (2) above implies that $d_E-r_Eg>i>0$ so
$\wedge^iM_E^\vee$ is well defined. From (1) the definition of $d_F$
makes sense and for $F\in \mathscr U(r_F,d_F)$,
$\mu(\wedge^iM_E^\vee \otimes F)=g-1$.

From \cite{gt}, Prop.~5.1, for a generic effective divisor $D_i$ of
degree $i$, there is an immersion ${\mathcal O}(D_i)\rightarrow
\wedge^iM_E^{\vee}$. Hence,
 it is enough to show that a general vector bundle $F\in \mathscr U(r_F,d_F)$ can be written as
 an extension of vector bundles
$ 0\to \mathscr O_C(-D_i) \to F \to \overline{F} \to 0$ .  Since
$i\ge g$,
 $\mathscr O_C(-D_i)$ is a general line bundle of degree $-i$. From
 \cite{rt}  a general $F\in\mathscr U(r_F,d_F)$ can be written as such an extension so long as
$ \mu(\overline{F})-\mu(\mathscr O_C(-D_i))\ge g-1 $. From the short
exact sequence above and the definition of $d_F$,
 $$\mu(\overline{F})-\mu(\mathscr O_C(-D_i))=\frac{r_F(g-1-i\frac{d_E}{d_E-r_Eg})+i}{ r_F-1}.$$
 Condition (2) ensures that the needed inequality is satisfied.
\end{proof}

\begin{cor}\label{cormu}
If $E\in \mathscr U(r_E,d_E)$ with $\mu(E)=\left(\frac{a}{b}\right)
g$ for some $a,b\in \mathbb N$, and $r_F,i\in \mathbb N$ satisfy the
following condtions:
\begin{enumerate}
\item  $r_Fi\frac{a}{a-b}\in \mathbb N$,

\item $(g-1)\frac{a}{b}\ge r_Fi\frac{a}{a-b}\ge r_F
 g\frac{a}{a-b}$,
 \end{enumerate}
then setting $d_F=r_F(g-1)-r_Fi\frac{a}{a-b}$,
$$\Theta_{\wedge^iM_E^\vee} :=\{F\in \mathscr U(r_F,d_F):h^0(F\otimes \wedge^iM_E^\vee)>0\}=\mathscr U(r_F,d_F).$$
     \end{cor}

\begin{rem}
With the  assumptions in the corollary,
 $\mu(\wedge^iM_E^\vee)=\frac{ia}{a-b} \equiv \frac{ib}{a-b} \pmod{\mathbb Z}$.
Set $r=\operatorname{rank}(\wedge^iM_E^\vee)$, and
$d=\deg(\wedge^iM_E^\vee)$. If follows that
$\frac{r}{\gcd(r,d)}=\frac{a-b}{\gcd(a-b,ia)}$. Set $ n= r_F\cdot
\frac{\gcd(a-b,ia)}{a-b}$. For any $L\in \operatorname{Pic}^d(C)$
and an appropriate line bundle $\eta$ on $C$, it follows that
$\wedge^i M_{E\otimes \eta}^\vee\in S\mathscr U(r,L)$ is a base
point for the linear system
 $|\mathscr L^{\otimes n}|$ on $S\mathscr U(r,L)$.  In other words, with the assumptions above
  $$\wedge^i M_{E\otimes \eta}^\vee\in \Theta_F
\subseteq S\mathscr U(r,L)\ \forall \ F\in \mathscr U(r_F,d_F).$$
\end{rem}

\begin{pro}\label{mainexa}
Fix $r_E,r_F,A,b,i'\in \mathbb N$ satisfying the following conditions
\begin{enumerate}
\item $b|r_Eg$
\item $\frac{g-1}{b}\ge i'\ge \frac{g}{A}$
\item $r_F\ge \frac{b}{A}$
\end{enumerate}
Set $n=\gcd(r_F,i'b)$, $r=\binom{\frac{Ag}{b}r_Fr_E}{i}$,
$d=r(i'A+\frac{i'b}{r_F})$,
 and fix a line bundle $L\in \operatorname{Pic}^d(C)$.
On $S\mathscr U(r,L)$ the dimension of the base locus of $|\mathscr L^{\otimes n}|$   is at least
$$
(r_E^2-1)(g-1)
$$
     \end{pro}
\begin{proof}Take  $a=Ar_F+b, i=i'A$ in the previous corollary.
Assumption (3) in the proposition insures that the map
 taking $E\mapsto M_E$ is generically finite (\cite{gt}).
\end{proof}

\begin{teo}\label{teo1}
Let $\mu\in \mathbb Q$, and suppose there exists numbers $i',r_F\in
\mathbb N$ with $1\le i'\le g-1$ such that $ \mu\equiv
\frac{i'}{r_F} \pmod{\mathbb Z}. $ Suppose that $N,r_E\in \mathbb
N$, $r_E\ge 2$, and $ N\ge N_0:=\binom{r_Fr_Eg^2}{i'g}$. For $L\in
\operatorname{Pic}^{\mu N}(C)$,  on $S\mathscr U (N,L)$, the base
locus of the linear system
 $|\mathscr L|$ has dimension at least
$ ((r_E^2-1)+(N-N_0)^2)(g-1)+1-\delta_{N,N_0}$. The same can be said
for the linear system $|\mathscr L^{\otimes
\gcd(r_F,i')}|$.\end{teo}

\begin{proof}
Take $A=g$ and $b=1$ in Proposition \ref{mainexa}.
\end{proof}

\begin{rem}
  For a given slope $\mu\in \mathbb Q$, taking a curve of high enough genus, the theorem implies that the base locus
   is nonempty provided the rank is sufficiently large.  In particular, Theorem \ref{teomain} follows
    from Theorem \ref{teo1}.
\end{rem}

Taking $A=2$, $b=1$, and $i'=r_F=g-1$ in Prop \ref{mainexa},   we
get examples with lower rank and higher level.

\begin{teo}\label{teo2} Suppose that $N,r_E\in \mathbb N$, $r_E\ge 2$, and
 $$ N\ge N_0:=\binom{2g(g-1)r_E}{2(g-1)}. $$ Then on $S\mathscr U (N,\mathscr O_C)$, the base locus of
  the linear system $|\mathscr L^{\otimes g-1}|$ has dimension at least $$ ((r_E^2-1)+(N-N_0)^2)(g-1)+1-\delta_{N,N_0}.
  $$\end{teo}


\end{document}